\documentclass{article}

\usepackage[utf8]{inputenc}

\usepackage{algorithm}
\usepackage{algorithmic}
\usepackage{amsmath}
\usepackage{amssymb}
\usepackage{amsthm}
\usepackage{arydshln}
\usepackage{authblk}

\usepackage{bbold}
\usepackage{bm}
\usepackage{booktabs}
\usepackage{bookmark}
\usepackage{calc}              
\usepackage{enumitem}
\usepackage{extarrows}
\usepackage{float}
\usepackage{geometry}
\usepackage{graphicx}          
\usepackage{hyperref}
\usepackage{indentfirst}
\usepackage{latexsym}
\usepackage{mathtools}
\usepackage{multirow}          
\usepackage{subcaption}
\usepackage{tikz}              
\usepackage{xcolor}            
\usetikzlibrary{shapes.geometric, arrows, chains}

\title{Robust space-time multiscale upscaling via multicontinuum homogenization for evolving perforated media}

\providecommand{\keywords}[1]{\textbf{Keywords:} #1}

\author[a,b]{Wei Xie}
\author[b]{Viet Ha Hoang}
\author[a]{Yin Yang\thanks{Corresponding author}}
\author[a]{Yunqing Huang}
\affil[a]{School of Mathematics and Computational Science, Xiangtan University, National Center for Applied Mathematics in Hunan, Xiangtan, Hunan, China 411105}
\affil[b]{Division of Mathematical Sciences, School of Physical and Mathematical Sciences, Nanyang Technological University, Singapore 637371}

\date{}

\begin{document}
\maketitle

\renewcommand{\thefootnote}{}
\footnotetext{E-mail addresses: xiew@smail.xtu.edu.cn (Wei Xie), 
vhhoang@ntu.edu.sg (Viet Ha Hoang),
yangyinxtu@xtu.edu.cn (Yin Yang),
huangyq@xtu.edu.cn (Yunqing Huang)}
\renewcommand{\thefootnote}{\arabic{footnote}}

\begin{abstract}
Time-evolving perforated domains arise in many engineering and geoscientific applications, including reactive transport, particle deposition, and structural degradation in porous media. Accurately capturing the macroscopic behavior of such systems poses significant computational challenges due to the dynamic fine-scale geometries. In this paper, we develop a robust and generalizable multiscale modeling framework based on multicontinuum homogenization to derive effective macroscopic equations in shrinking domains. The method distinguishes multiple continua according to the physical characteristics (e.g., channel widths), and couples them via space-time local cell problems formulated on representative volume elements. These local problems incorporate temporal derivatives and domain evolution, ensuring consistency with underlying fine-scale dynamics. The resulting upscaled system yields computable macroscopic coefficients and is suitable for large-scale simulations. Several numerical experiments are presented to validate the accuracy, efficiency, and potential applicability of the method to complex time-dependent engineering problems.
\end{abstract}

\section{Introduction}
Many natural processes involve complex interactions across multiple spatial and temporal scales. For instance, mineral dissolution and precipitation in porous media progressively alter pore structures, thereby affecting fluid flow patterns~\cite{devigne2008numerical, bringedal2020phase}. Similarly, the injection of CO$_2$ or water may induce mechanical failure in initially intact solids, leading to crack formation and structural deformation~\cite{miehe2016phaseIII}. Particle deposition is another example, in which pores become gradually clogged, dynamically reshaping the flow domain~\cite{wang2018numerical}. A common characteristic of these processes is that the computational domain evolves over time, posing significant challenges for accurate and efficient numerical simulation. Capturing such dynamics requires advanced multiscale methodologies that couple microscale evolution with macroscale responses.

A widely adopted approach for multiscale modeling involves two-level discretization strategies. These methods construct localized basis functions on fine grids to resolve subgrid-scale features while solving a coarse-scale problem on a reduced mesh to improve computational efficiency. Representative techniques include the multiscale finite element method (MsFEM)~\cite{hou1997multiscale, le2014msfem}, the generalized multiscale finite element method (GMsFEM)~\cite{efendiev2013generalized, alikhanov2025multiscale, chung2018generalized, chung2018space, xie2025multiscale}, wavelet-based edge multiscale finite element methods~\cite{fu2019edge, fu2021wavelet}, the constraint energy minimizing GMsFEM (CEM-GMsFEM)~\cite{chung2018constraint, ye2023constraint, chung2021convergence, xie2024cem}, and the nonlocal multicontinuum method (NLMC)~\cite{chung2018non, leung2019space, hu2025space}. These approaches allow for accurate resolution of fine-scale heterogeneities while maintaining global solution quality. For space-time multiscale problems, these methods have been extended to construct space-time basis functions~\cite{chung2018generalized, chung2018space, leung2019space, hu2025space}.

Another class of methods seeks to derive effective macroscopic models by computing homogenized coefficients that capture microscale behavior. Representative approaches include numerical homogenization~\cite{jikov2012homogenization, pavliotis2008multiscale, efendiev2004numerical, altmann2021numerical} and the heterogeneous multiscale method~\cite{abdulle2012heterogeneous, ming2005analysis}. These methods typically rely on the assumption of scale separation in space and/or time, and often make use of asymptotic expansions, such as the two-scale ansatz.

Recent advances in multicontinuum homogenization~\cite{efendiev2023multicontinuum, efendiev2024multicontinuum, chung2024multicontinuum, leung2024some, bai2025multicontinuum, xie2025multicontinuum, xie2025hierarchical} have extended the applicability of homogenization to settings without clear scale separation. Instead of asymptotic expansions, these approaches use multiple localized basis functions and smooth macroscopic variables to describe the dynamics of distinct continua within a representative volume element (RVE). Continua are defined based on physical attributes such as material contrast or channel width, and their interactions are modeled via carefully formulated local constraint problems.

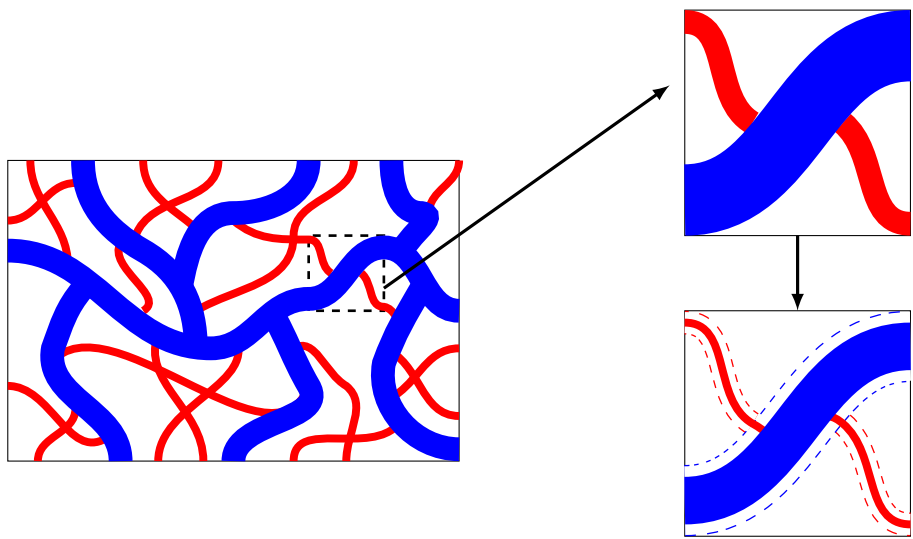
\begin{figure}[H]
\centering
\begin{tikzpicture}
\draw (0,0) rectangle (6,4);
\draw[dashed,line width=1 pt] (4,2) rectangle (5,3);
\draw (9,3) rectangle (12,6);
\draw (9,-1) rectangle (12,2);
\pgfmathsetmacro{\a}{9};
\pgfmathsetmacro{\b}{\a*0.0351/2};
\pgfmathsetmacro{\c}{\a/3};
\pgfmathsetmacro{\d}{\b/3};
\draw[color=red,line width=\c pt] (4.3,2.5) to[out=150,in=0]  (4,3-\d)
to[out=180,in=320] (3,3.2);
\draw[color=red,line width=\c pt] (4.7,2.5) to[out=320,in=180]  (5,2+\d)
to[out=350,in=160] (5.3,1.8);
\draw[color=red,line width=\c pt] (0,1) to[out=0,in=210]  (1.3,0.3);
\draw[color=red,line width=\c pt] (0.4,0) to[out=90,in=270]  (0.9,0.7);
\draw[color=red,line width=\c pt] (0,3.2) to[out=0,in=180]  (0.9,3.7);
\draw[color=red,line width=\c pt] (0.8,2.6) to[out=90,in=270]  (0.3,4);
\draw[color=red,line width=\c pt] (1.5,2.9) to[out=240,in=30]  (1.8,2);
\draw[color=red,line width=\c pt] (2.5,3.1) to[out=140,in=270]  (1.8,4);
\draw[color=red,line width=\c pt] (1.7,3) to[out=70,in=270]  (2.8,4);
\draw[color=red,line width=\c pt] (5.6,3.3) to[out=90,in=270]  (6,4);
\draw[color=red,line width=\c pt] (3.8,1.5) to[out=30,in=150] (4.5,1)
to [out=320,in=90]  (4.5,0);
\draw[color=red,line width=\c pt] (5.2,0.5) to[out=230,in=90]  (3.8,0);
\draw[color=red,line width=\c pt] (2.5,2) to[out=30,in=270]  (3.8,3)
to[out=90,in=270] (4.6,4);
\draw[color=red,line width=\c pt] (5,0.8) to[out=20,in=180]  (6,1.5);
\draw[color=red,line width=\c pt] (5,1.5) to[out=340,in=180]  (6,0.6);
\draw[color=red,line width=\c pt] (0.5,1.2) to[out=60,in=250]  (4,1);
\draw[color=red,line width=\c pt] (2,1.7) to[out=250,in=90]  (2.6,0);
\draw[color=red,line width=\c pt] (3,1.7) to[out=300,in=90]  (2,0);
\draw[color=red,line width=3*\c pt] (9.9,4.5) to[out=150,in=0]  (9,6-3*\d);
\draw[color=red,line width=3*\c pt] (11.1,4.5) to[out=320,in=180]  (12,3+3*\d);
\draw[dashed,color=red,line width=3*\c pt] (9.9,0.5) to[out=150,in=0]  (9,2-3*\d);
\draw[dashed,color=red,line width=3*\c pt] (11.1,0.5) to[out=320,in=180]  (12,-1+3*\d);
\draw[color=white,line width=3*\c-1 pt] (9.9,0.5) to[out=150,in=0]  (9,2-3*\d);
\draw[color=white,line width=3*\c-1 pt] (11.1,0.5) to[out=320,in=180]  (12,-1+3*\d);
\draw[color=blue,line width=\a pt] (0,2.8) to[out=0,in=180]  (2.7,1.5) 
to[out=0,in=180] (4,2+\b) to[out=0,in=180] (5,3-\b)
to[out=0,in=180] (6,2);
\draw[color=blue,line width=\a pt] (1.5,0) to[out=90,in=260]  (0.6,1.5)
to[out=70,in=230] (1,2.3);
\draw[color=blue,line width=\a pt] (1,4) to[out=270,in=150]  (1.8,2.8)
to[out=330,in=90] (2.5,1.5);
\draw[color=blue,line width=\a pt] (3,0) to[out=90,in=300]  (4,1)
to[out=120,in=300] (3.5,2);
\draw[color=blue,line width=\a pt] (4,4) to[out=270,in=0]  (3,3.3)
to[out=180,in=80] (2.3,2.3);
\draw[color=blue,line width=\a pt] (6,\b) to[out=180,in=260]  (5,1.3)
to[out=70,in=240] (5.5,2.3);
\draw[color=blue,line width=\a pt] (5.1,4) to[out=270,in=180]  (5.5,3.3)
to[out=350,in=80] (5.2,2.7);
\draw[color=blue,line width=3*\a pt] (9,3+3*\b) to[out=0,in=180]  (12,6-3*\b);
\draw[dashed,color=blue,line width=3*\a pt] (9,-1+3*\b) to[out=0,in=180]  (12,2-3*\b);
\draw[color=white,line width=3*\a-1 pt] (9,-1+3*\b) to[out=0,in=180]  (12,2-3*\b);
\draw[color=red,line width=\c pt] (10.2,0.1) to[out=80,in=330] (9.9,0.5) 
to[out=150,in=0]  (9,2-3*\d);
\draw[color=red,line width=\c pt] (10.9,0.6) to[out=350,in=140] (11.1,0.5) to[out=320,in=180]  (12,-1+3*\d);
\draw[color=blue,line width=2*\a pt] (9,-1+3*\b) to[out=0,in=180]  (12,2-3*\b);
\draw[very thick,-latex] (5,2.3) -- (8.8,5);
\draw[very thick,-latex] (10.5,3) -- (10.5,2);
\end{tikzpicture}
\caption{Illustration of shrinking domain.}
\label{fig:shrinkingriver}
\end{figure}

In this work, we present a novel space-time multicontinuum homogenization framework tailored for time-dependent problems in shrinking perforated domains. These domains evolve due to processes such as pore clogging or crack propagation. As a representative example, we consider a two-continuum system, though the proposed framework naturally extends to more general multichannel configurations. As illustrated in Figure~\ref{fig:shrinkingriver}, the blue region denotes the first continuum (thick channels), and the red region indicates the second continuum (thin channels). Each continuum $\psi_i(\boldsymbol{x},t)$ is represented by a space-time function that equals one in the $i$th continuum and zero elsewhere. To resolve fine-scale behavior, we solve local constraint problems on oversampled RVEs, which mitigate boundary artifacts and enable accurate computation of effective properties.

The main contributions of this paper are as follows: We introduce a space-time multicontinuum homogenization framework tailored for evolving perforated domains. This framework involves the construction of local cell problems that account for both temporal derivatives and geometric evolution. From the local solutions, we derive computable macroscopic equations that describe coarse-scale dynamics. Finally, we demonstrate the accuracy and efficiency of the proposed method through extensive numerical experiments.

The rest of the paper is organized as follows. Section~\ref{sec:preliminaries} introduces the model problem and fine-scale discretization. Section~\ref{sec:mh} presents the multicontinuum homogenization framework. Numerical results are reported in Section~\ref{sec:numuerical_ex}, and conclusions are provided in Section~\ref{sec:conclusions}.

\section{Preliminaries} \label{sec:preliminaries}

In this section, we introduce the model problem, establish the relevant notations, and describe the fine-scale discretization.

We consider the following parabolic equation posed in a time-evolving perforated domain:
\begin{equation}
\begin{cases}
\begin{aligned}
\frac{\partial u(\boldsymbol{x},t)}{\partial t} - \mathrm{div} \left( \kappa(\boldsymbol{x}) \nabla u(\boldsymbol{x},t) \right) &= f(\boldsymbol{x},t), 
& \text{in}~ \Omega^{\epsilon}(t) \times I, \\
u(\boldsymbol{x},0) &= u_0(\boldsymbol{x}), 
& \text{in}~ \Omega^{\epsilon}(0), \\
u(\boldsymbol{x},t) &= 0, 
& \text{on}~ \partial\Omega^{\epsilon}(t) \times I,
\end{aligned}
\end{cases}
\label{eq:pde}
\end{equation}
where $\Omega^{\epsilon}(t)$ denotes the computational domain at time $t$, evolving due to physical processes such as pore clogging or particle deposition. The time interval is denoted by $I = (0, T]$. The parameter $\epsilon$ characterizes the microscale heterogeneity of the domain, such as the smallest pore or channel width. Although $\epsilon$ may depend on time, we omit this dependency for simplicity of notation.

Let $H^1(\Omega^{\epsilon})$ denote the standard Sobolev space. We define the subspace
\[
H_0^1(\Omega^{\epsilon}(t)) \coloneqq 
\left\{ v \in H^1(\Omega^{\epsilon}(t)) \mid v = 0 \text{ on } \partial \Omega^{\epsilon}(t) \right\},
\]
and introduce the bilinear forms
\[
a(u,v) = \int_{\Omega^{\epsilon}} \kappa \nabla u \cdot \nabla v, \quad
(w,v) = \int_{\Omega^{\epsilon}} w v.
\]

The weak formulation of Equation~\eqref{eq:pde} reads: find $u \in H_0^1(\Omega^{\epsilon}(t))$ such that
\begin{align}
\left( \frac{\partial u}{\partial t}, v \right) + a(u, v) &= (f, v), \quad \forall v \in H_0^1(\Omega^{\epsilon}(t)), 
\label{eq:weak_pde} \\
(u, v) &= (u_0, v), \quad \forall v \in H_0^1(\Omega^{\epsilon}(t)).
\label{eq:weak_init}
\end{align}

We discretize the time interval using $0 = t_0 < t_1 < \cdots < t_N = T$, and denote the time step by $\tau_k = t_{k+1} - t_k$, with $\tau = \min_{1 \le k \le N} \tau_k$. We assume that the time step is sufficiently small to resolve the dynamic evolution of the domain.

Let $\mathcal{T}_h$ be a conforming partition of the initial domain $\Omega^{\epsilon}(0)$ into triangles or quadrilaterals. The mesh size $h$ is chosen to be fine enough to resolve both the microscale coefficients and the geometric complexity of the domain, thereby avoiding the failure of standard finite element methods.

We define $V_h$ as the space of continuous, piecewise linear functions associated with the mesh $\mathcal{T}_h$. Due to the time-dependent nature of the computational domain, not all basis functions in $V_h$ remain active at each time step. For each $t_k$, we define the active finite element space $V_h^k$ as the subspace of $V_h$ spanned by basis functions whose associated nodes lie in the interior of $\Omega^{\epsilon}(t_k)$. In other words, the degrees of freedom of the solution $u^k$ are restricted to the current domain $\Omega^{\epsilon}(t_k)$, reflecting the evolution of the perforated geometry.

To discretize Equation~\eqref{eq:weak_pde} in time, we employ the backward Euler method. For $k=0,\ldots,N-1$, we seek $u^{k+1} \in V_h^{k+1}$ such that
\begin{equation}
\left( \frac{u^{k+1} - u^k}{\tau_k}, v \right) + a(u^{k+1}, v) = (f, v), 
\quad \forall v \in V_h^{k+1}.
\end{equation}

The initial condition is enforced by projecting $u_0$ onto the discrete space defined on the initial domain. That is, we compute $u_h^0 \in V_h^0$ such that
\begin{equation}
(u_h^0, v) = (u_0, v), \quad \forall v \in V_h^0.
\end{equation}
This provides a consistent initialization of the fine-grid solution at $t = 0$.

\section{Multicontinuum Homogenization} \label{sec:mh}

\begin{figure}[H]
\centering
\begin{tikzpicture}
\node[anchor=south west,inner sep=0] at (9.2,2.2) {\includegraphics[width=0.6cm]{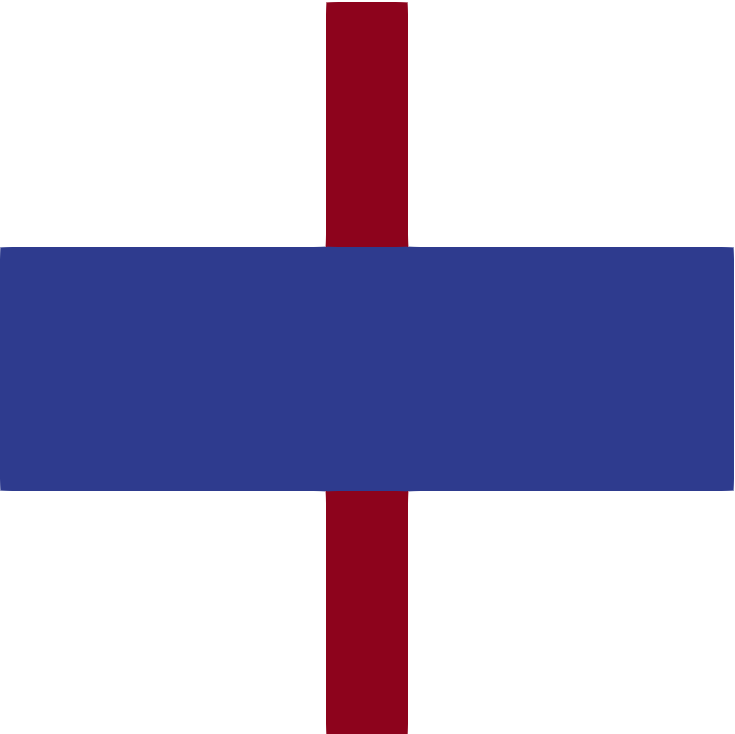}};
\draw[help lines,step=0.5cm] (0,0) grid (6,5);
\draw[ultra thick,->] (4.5,1.25) -- (7.8,2.5);
\draw[->] (9.8,2.3) -- (11.8,1.5);
\draw[->] (10.5,2.9) -- (11.8,2.5);
\draw[->] (11,3.5) -- (11.8,3.5);
\draw[ultra thick] (4,1) rectangle (4.5,1.5);
\draw (8,1) rectangle (11,4);
\draw (8.5,1.5) rectangle (10.5,3.5);
\draw (9.2,2.2) rectangle (9.8,2.8);
\node[anchor=south] at (6.5,0) {$\Omega^{\epsilon}$};
\node at (12.2,3.5) {$K_p^{\epsilon}$};
\node at (12.2,2.5) {$R_p^+$};
\node at (12.2,1.5) {$R_p$};
\end{tikzpicture}
\caption{Illustration of the representative volume element (RVE) $R_p$, oversampling domain $R_p^+$, and coarse block $K_p^{\epsilon}$.}
\label{fig:macropic_point_rve}
\end{figure}
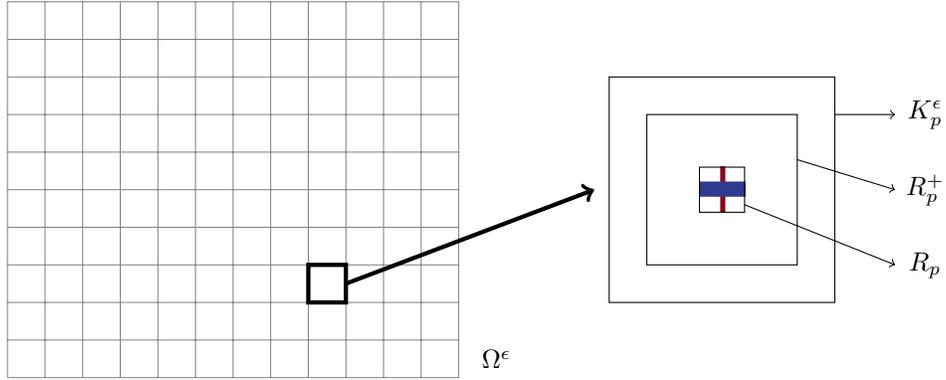

In this section, we employ the multicontinuum homogenization framework to derive the macroscopic equations associated with the original problem~\eqref{eq:pde}.

We begin by introducing the computational setting for the multicontinuum homogenization procedure. Let $\mathcal{T}_H$ denote a coarse partition of the global domain $\Omega^{\epsilon}$, where $H$ is the coarse grid size, assumed to be larger than the characteristic scale of heterogeneity $\epsilon$ (e.g., the typical pore or perforation width). 
For each coarse block $K_p^{\epsilon} \in \mathcal{T}_H$, we define a representative volume element (RVE) $R_p \subset K_p^{\epsilon}$. 
To mitigate artificial boundary effects~\cite{hou1997multiscale,chung2018constraint,xie2024cem}, we construct an oversampled region $R_p^+ \supset R_p$, within which the local problems are solved. Figure~\ref{fig:macropic_point_rve} illustrates the relationship between $K_p^{\epsilon}$, $R_p$, and $R_p^+$.
We distinguish several scales: $h < \epsilon < \varepsilon < H$, where $h$ is the fine-grid size, $\epsilon$ is the heterogeneity scale, $\varepsilon$ is the RVE scale, and $H$ is the coarse-grid size.

We now present a multiscale expansion of the fine-scale solution $u$ within each representative volume element (RVE). 
Following our previous work~\cite{xie2025multicontinuum}, we distinguish multiple continua based on their geometric or physical characteristics (e.g., channel width); see Figure~\ref{fig:shrinkingriver} for an illustration.
The solution $u$ is approximated by a combination of smooth macroscopic variables $U$, which represent large-scale trends, and localized basis functions $\phi$, which are designed to capture fine-scale features while preserving coarse-scale averages.
Specifically, we consider the expansion
\begin{equation}
u = \phi_i U_i + \phi_i^m \frac{\partial U_i}{\partial x_m}
+ \phi_i^{mn} \frac{\partial^2 U_i}{\partial x_m \partial x_n} + \cdots,
\label{eq:expansion_u_rve}
\end{equation}
where $U_i$ denotes the macroscopic variable associated with the $i$-th continuum, and a set of localized basis functions $\phi$ are constructed to preserve the average behavior of $u$ while encoding fine-scale effects. 
The concepts behind the macroscopic variables $U$ and the definitions of the basis functions $\phi$ will be discussed later.

In this work, we retain only the first two terms in the expansion~\eqref{eq:expansion_u_rve}, leading to the approximation:
\begin{equation}
u \approx \phi_i U_i + \phi_i^m \frac{\partial U_i}{\partial x_m}.
\label{eq:u_expansion}
\end{equation}

In the context of the multicontinuum homogenization framework~\cite{efendiev2023multicontinuum,xie2025multicontinuum}, the macroscopic variable $U_i$ is defined as the weighted average of the fine-scale solution $u$ over the $i$-th continuum. 
Given an RVE $R_p$, this average is expressed as
\[
U_i(\boldsymbol{x}_p^*) = 
\frac{\int_{R_p} u \psi_i}{\int_{R_p} \psi_i},
\]
where $\boldsymbol{x}_p^*$ denotes the centroid of $R_p$, and $\psi_i(\boldsymbol{x}, t)$ is a characteristic function that equals one inside the $i$-th continuum and zero elsewhere.

The local basis functions $\phi_i$ and $\phi_i^m$ are obtained by solving constrained local problems that enforce consistency with constant and linear macroscopic fields. These problems ensure that the local averages of the basis functions match the intended macroscopic behavior.

To capture the constant behavior in each continuum, the basis functions $\phi_i$ are computed by solving:
\begin{equation}
\begin{cases}
\begin{aligned}
\int_{R_p^+} \frac{\partial \phi_i}{\partial t} v +
\int_{R_p^+} \kappa \nabla \phi_i \cdot \nabla v -
\sum_{j,q} \frac{\beta_{ij}^q}{\int_{R_q} \psi_j}
\int_{R_q} \psi_j v = 0, &
&\text{in}~ (0,T], \\
\int_{R_p^+} \kappa \nabla \phi_i \cdot \nabla v -
\sum_{j,q} \frac{\zeta_{ij}^q}{\int_{R_q} \psi_j}
\int_{R_q} \psi_j v = 0, &
&\text{at}~ t=0, \\
\int_{R_q} \phi_i \psi_j
= \delta_{ij} \int_{R_q} \psi_j,&
&\text{for all}~ t \in [0,T], \\
\end{aligned}
\end{cases}
\label{eq:phi_avar}
\end{equation}
where $\beta$ and $\zeta$ are Lagrange multipliers enforcing the averaging constraints, ensuring that each $\phi_i$ has unit mean over the $i$-th continuum and zero mean elsewhere.

To capture linear variations in the macroscale, the basis functions $\phi_i^m$ are computed by solving:
\begin{equation}
\begin{cases}
\begin{aligned}
\int_{R_p^+} \frac{\partial \phi_i^m}{\partial t} v +
\int_{R_p^+} \kappa \nabla \phi_i^m \cdot \nabla v -
\sum_{j,q} \frac{\beta_{ij}^{mq}}{\int_{R_q} \psi_j}
\int_{R_q} \psi_j v = 0,&
&\text{in}~ (0,T], \\
\int_{R_p^+} \kappa \nabla \phi_i^m \cdot \nabla v -
\sum_{j,q} \frac{\zeta_{ij}^{mq}}{\int_{R_q} \psi_j}
\int_{R_q} \psi_j v = 0,&
&\text{at}~ t=0, \\
\int_{R_q} \phi_i^m \psi_j = \delta_{ij} \int_{R_q} (x_m-c_{mj})\psi_j,&
&\text{for all}~ t \in [0,T]. \\
\end{aligned}
\end{cases}
\label{eq:phi_grad}
\end{equation}
Here, $c_{mj}$ is chosen such that $\int_{R_q} (x_m - c_{mj}) \psi_j = 0$, ensuring zero average deviation within each continuum.
Again, the Lagrange multipliers enforce the constraints so that $\phi_i^m$ matches the mean of $x_m$ in the $i$-th continuum and vanishes in the others.

The resulting basis functions satisfy the following scaling properties (see~\cite{efendiev2023multicontinuum} for details):
\begin{equation}
\begin{aligned}
\Vert \phi_i \Vert = \mathcal{O}(1), \quad
\Vert \nabla \phi_i \Vert = \mathcal{O}\left(\frac{1}{\varepsilon}\right), \\
\Vert \phi_i^m \Vert = \mathcal{O}(\varepsilon), \quad
\Vert \nabla \phi_i^m \Vert = \mathcal{O}(1).
\end{aligned}
\label{eq:phi_scaling}
\end{equation}

We now derive the macroscopic equations by substituting the multiscale expansion~\eqref{eq:u_expansion} into the weak formulation~\eqref{eq:weak_pde}.
Summing over all RVEs and applying a standard localization argument yields:
\begin{equation}
\sum_p \frac{|K_p^{\epsilon}|}{|R_p|} \int_{R_p} 
\frac{\partial u}{\partial t} v 
+ \sum_p \frac{|K_p^{\epsilon}|}{|R_p|} \int_{R_p} 
\kappa \nabla u \cdot \nabla v
= \sum_p \frac{|K_p^{\epsilon}|}{|R_p|} \int_{R_p} 
f v.
\label{eq:sum_integ_rve}
\end{equation}

Using the multiscale approximation~\eqref{eq:u_expansion} and the assumption that $U_i$ varies slowly within $R_p$ compared to the local basis functions $\phi$, we obtain:
\begin{equation}
\int_{R_p} \frac{\partial (\phi_i U_i)}{\partial t} v
+ \int_{R_p} \frac{\partial_m (\phi_i^m U_i)}{\partial t} v
\approx U_i(x_p^*) \int_{R_p} \frac{\partial \phi_i}{\partial t} v
+ \partial_m U_i(x_p^*) \int_{R_p} \frac{\partial \phi_i^m}{\partial t} v,
\label{eq:phiU_tapprox}
\end{equation}
and
\begin{equation}
\int_{R_p} \kappa \nabla (\phi_i U_i) \cdot \nabla v 
+ \int_{R_p} \kappa \nabla (\phi_i^m \frac{\partial U_i}{\partial x_m}) \cdot \nabla v 
\approx U_i(x_p^*) \int_{R_p} \kappa \nabla \phi_i \cdot \nabla v 
+ \frac{\partial U_i}{\partial x_m}(x_p^*) \int_{R_p} \kappa \nabla \phi_i^m \cdot \nabla v.
\label{eq:phiU_kapprox}
\end{equation}

By using the same approximation for the test function $v \approx \phi_j V_j + \phi_j^n \frac{\partial V_j}{\partial x_n}$, and applying~\eqref{eq:phiU_tapprox}--\eqref{eq:phiU_kapprox}, we obtain the upscaled variational problem, 

\[
\begin{aligned}
\sum_p \frac{|K_p^{\epsilon}|}{|R_p|} 
\left(V_j \frac{\partial U_i}{\partial t} 
\int_{R_p} \phi_i \phi_j
+  \frac{\partial U_i}{\partial_t} \frac{\partial V_j}{\partial x_n}
\int_{R_p} \phi_i \phi_j^n
+ V_j \frac{\partial U_i}{\partial_t \partial x_m}
\int_{R_p} \phi_i^m \phi_j 
+ \frac{\partial V_j}{\partial x_n} \frac{\partial U_i}{\partial_t \partial x_m}
\int_{R_p} \phi_i^m \phi_j^n
\right. \\
\left.
+ U_i V_j \int_{R_p} \phi_i \phi_j
+ U_i \frac{\partial V_j}{\partial x_n} \int_{R_p} \kappa \nabla \phi_i \phi_j^n
+ V_j \frac{\partial U_i}{\partial x_m} \int_{R_p} \phi_i^m \phi_j 
+ \frac{\partial U_i}{\partial x_m} \frac{\partial V_j}{\partial x_n} \int_{R_p} \phi_i^m \phi_j^n 
\right) \\
= \sum_p \frac{|K_p^{\epsilon}|}{|R_p|} \left(
V_j \int_{R_p} f \phi_j + \partial_n V_j \int_{R_p} f \phi_j^n
\right).
\end{aligned}
\]
which leads to the following macroscopic system in strong form:
\begin{equation}
D_{ji} \frac{\partial U_i}{\partial t} +
B_{ji} U_i 
- \partial_n (B_{ji}^{mn} \partial_m U_i) = b_j,
\label{eq:macro_eq_pde}
\end{equation}
Here, other terms will be ignored by the symmetric properties or the scaling in the \eqref{eq:phi_scaling}.
The effective coefficients given by
\[
D_{ji} = \int_{R_p} \phi_i \phi_j, \quad
B_{ji} = \int_{R_p} \kappa \nabla \phi_i \cdot \nabla \phi_j, \quad
B_{ji}^{mn} = \int_{R_p} \kappa \nabla \phi_i^m \cdot \nabla \phi_j^n, \quad
b_j = \int_{R_p} f \phi_j.
\]

Follow the same procedure, we obtain the coarse scale upscaled equation for initial condition,
\begin{equation}
D_{ji} U_i 
- \partial_n(D_{ji}^{mn} \partial_m U_i) = b_j,
\label{eq:macro_eq_init}
\end{equation}
where
\[
D_{ji} = \int_{R_p} \phi_i \phi_j, \quad
D_{ji}^{m} = \int_{R_p} \phi_i^m \phi_j, \quad
b_j = \int_{R_p} u_0 \phi_j.
\]

By combining~\eqref{eq:macro_eq_pde} and~\eqref{eq:macro_eq_init}, we obtain the complete macroscopic model for the multicontinuum system. It is important to note that the upscaled model is defined on the domain $\Omega$.

\section{Numerical Experiments} \label{sec:numuerical_ex}

In this section, we present three numerical experiments to demonstrate the accuracy and efficiency of the proposed multiscale method. The computational domain $\Omega$ is taken as the unit square. Each coarse block is treated as a RVE, and the fine grid resolution is set to $h = 1/1260$. We consider two coarse grid sizes, namely $H = 1/10$ and $H = 1/20$. For the oversampling strategy, we use 2 coarse layers when $H = 1/10$ and 4 layers when $H = 1/20$.

The domain shrinkage is modeled by prescribing a fixed contraction rate for each continuum. At every time step, the first continuum (corresponding to thick channels) contracts by three fine grid cells, whereas the second continuum (representing thin channels) contracts by one fine grid cell.

To assess the performance of the proposed method, we compute the relative $L^2$ error for each continuum using the following formula:
\begin{equation}
e_2^{(i)} = \frac{\sum_p \left| \frac{1}{|K_p|} \int_{K_p} U_i - \frac{1}{|{K_p^{\epsilon}} \cap \Omega_i|} \int_{{K_p^{\epsilon}}\cap \Omega_i} u \right|^2}{\sum_p \left| \frac{1}{|{K_p^{\epsilon}} \cap \Omega_i|} \int_{{K_p^{\epsilon}}\cap \Omega_i} u \right|^2},
\end{equation}
where $U_i$ denotes the multiscale solution for continuum $i$, and $u$ is the corresponding fine-scale reference solution.

The permeability field is defined as $\kappa = \kappa_1 \kappa_2$, where $\kappa_1$ characterizes the contrast between the continua: $\kappa_1 = 1$ in the first continuum and $\kappa_1 = 10^{-2}$ in the second. The component $\kappa_2$ varies across the three examples as follows:
\begin{itemize}
  \item \textbf{Example 1}: $\kappa_2 = 1$;
  \item \textbf{Example 2}: $\kappa_2 = 2 + \sin(2\pi x_1)\sin(20\pi x_2)$;
  \item \textbf{Example 3}: $\kappa_2 = 2 + \sin(\sqrt{20} \pi x_1)\sin(\pi x_2)$.
\end{itemize}

The source term is defined as
\[
f(x,y) = 10^{-3} \exp\left(-100\left((x - 0.5)^2 + (y - 0.5)^2\right)\right),
\]
and the initial condition is given by
\[
u_0(x,y) = 10^{-1} \exp\left(-100\left((x - 0.5)^2 + (y - 0.5)^2\right)\right).
\]

\subsection{Example 1}

Table~\ref{tab:ers_ex1} reports the numerical errors for various coarse grid sizes $H$. As expected, the error consistently decreases with mesh refinement, confirming the convergence of the proposed multiscale method. Moreover, the method maintains stable accuracy over the full simulation time horizon.

\begin{figure}[H]
\centering
\begin{subfigure}[b]{0.45\linewidth}
\includegraphics[width=\linewidth]{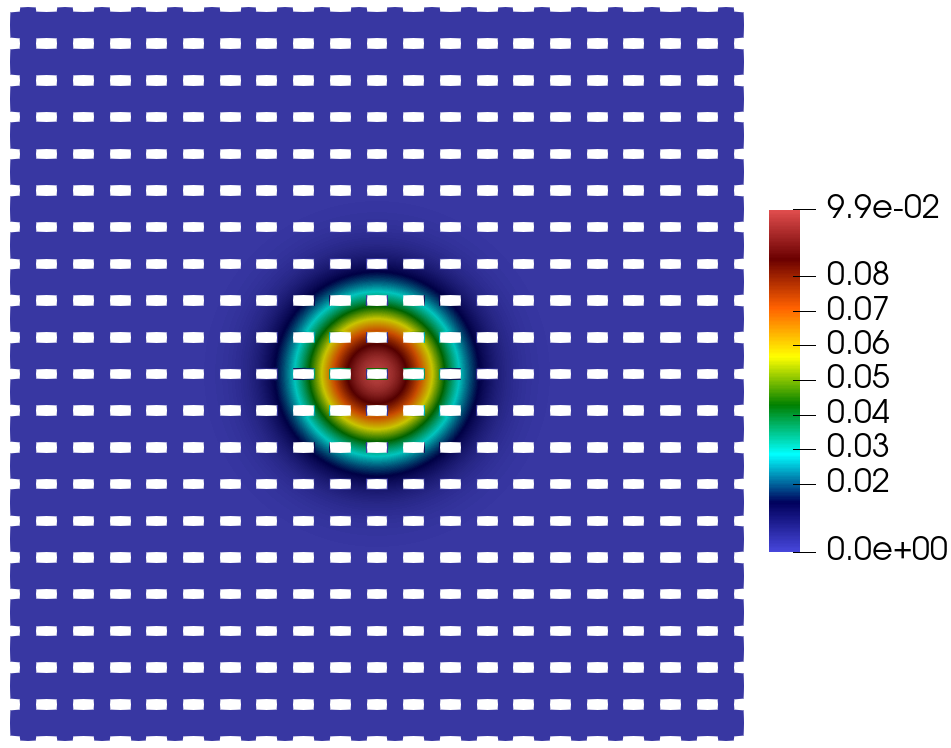}
\caption{$u_h^0$}
\end{subfigure}
\hfill
\begin{subfigure}[b]{0.45\linewidth}
\includegraphics[width=\linewidth]{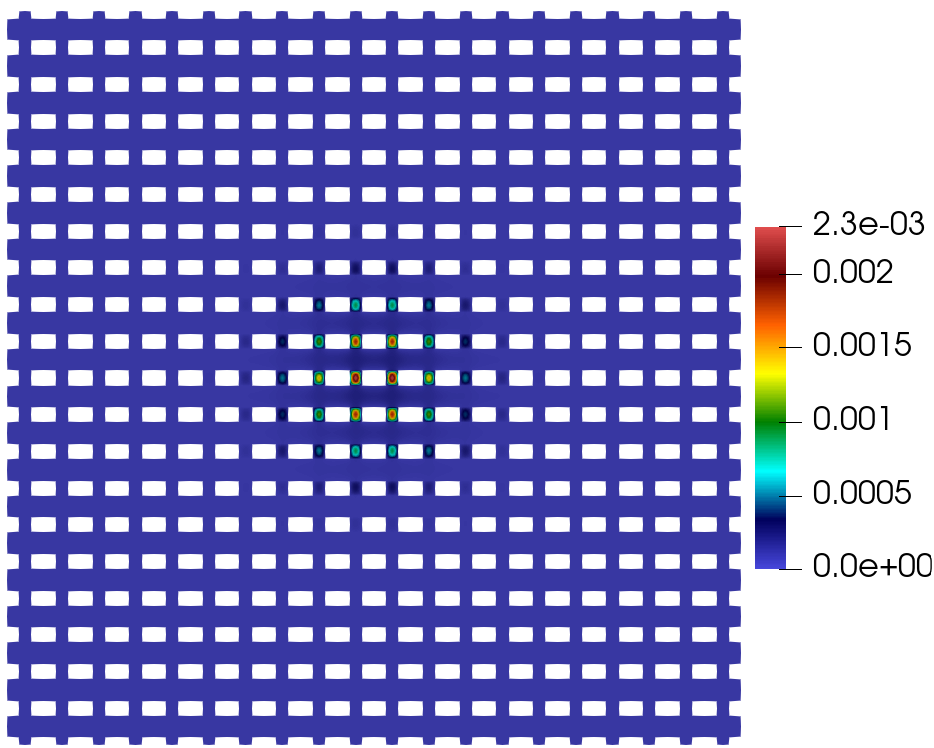}
\caption{$u_h^1$}
\end{subfigure}
\\[0.5cm]
\begin{subfigure}[b]{0.45\linewidth}
\includegraphics[width=\linewidth]{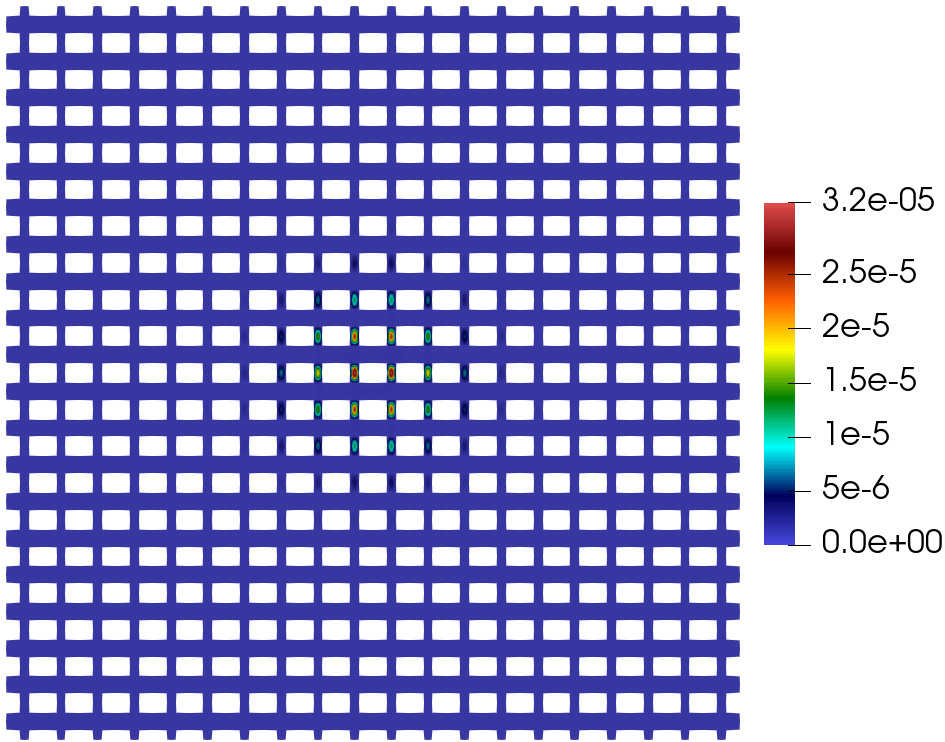}
\caption{$u_h^2$}
\end{subfigure}
\hfill
\begin{subfigure}[b]{0.45\linewidth}
\includegraphics[width=\linewidth]{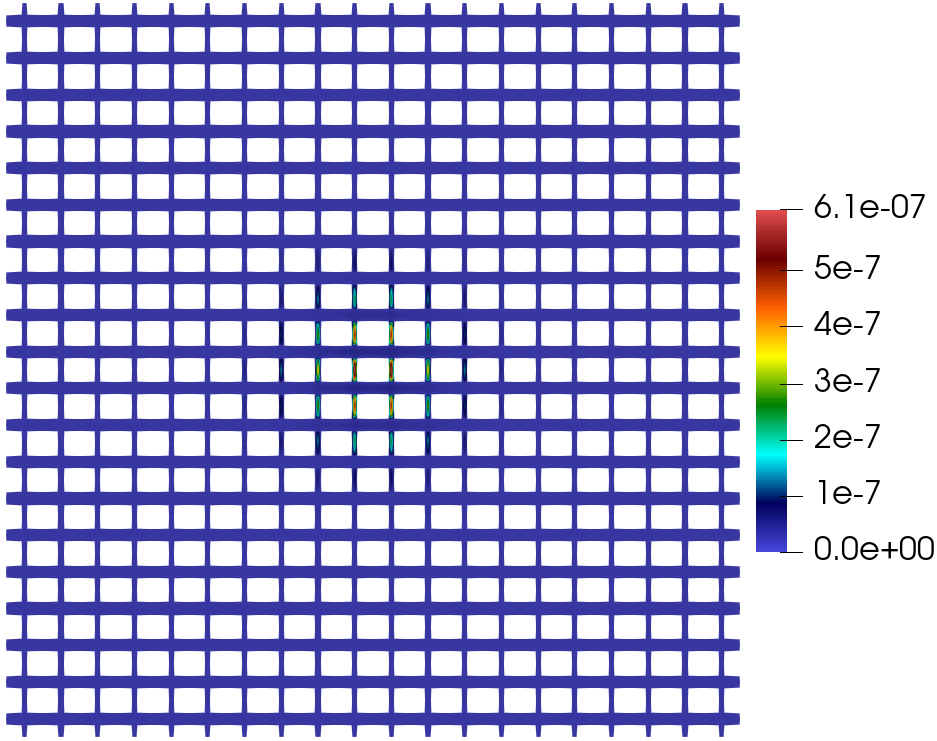}
\caption{$u_h^3$}
\end{subfigure}
\caption{Fine-grid solutions at different time steps in Example 1.}
\label{fig:ex1_uh}
\end{figure}

Figure~\ref{fig:ex1_uh} presents snapshots of the reference fine-scale solution at selected time instances. The solution exhibits higher magnitude and longer residence time within the thin channels, where the permeability is relatively low, resulting in slower fluid transport. In contrast, flow within the thick high-permeability channels is more transient, as the fluid quickly propagates and redistributes into adjacent regions. This spatial distribution reflects the influence of the heterogeneous permeability structure on the transport dynamics.

\begin{figure}[H]
\centering
\begin{subfigure}[b]{0.45\linewidth}
\includegraphics[width=\linewidth]{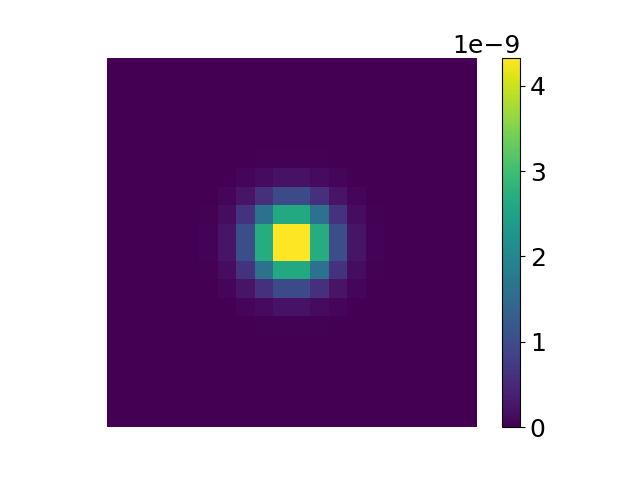}
\caption{Fine-grid average in continuum 1.}
\end{subfigure}
\hfill
\begin{subfigure}[b]{0.45\linewidth}
\includegraphics[width=\linewidth]{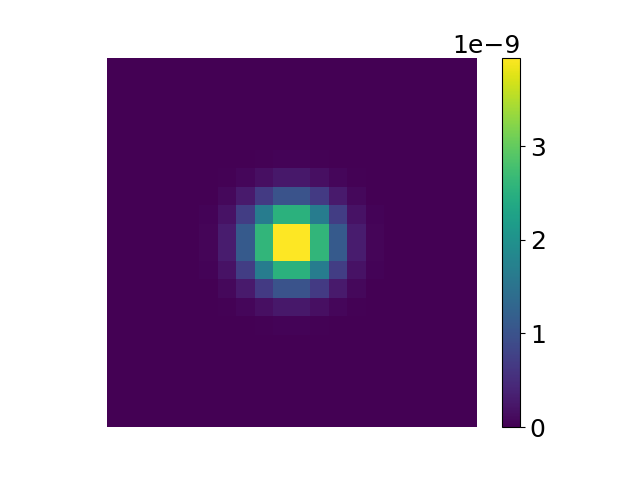}
\caption{Multiscale average in continuum 1.}
\end{subfigure}
\\[0.5cm]
\begin{subfigure}[b]{0.45\linewidth}
\includegraphics[width=\linewidth]{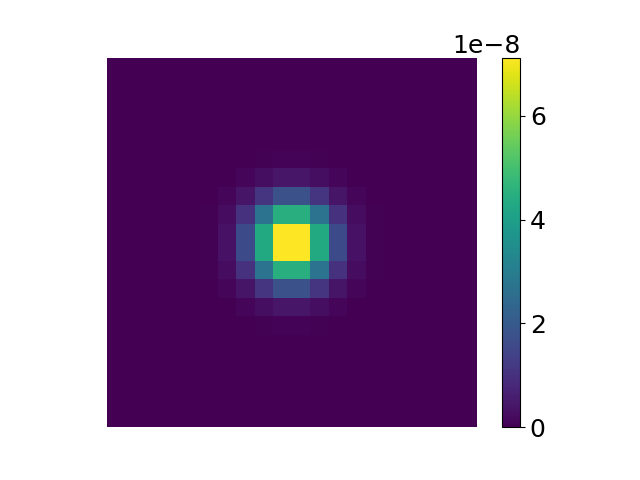}
\caption{Fine-grid average in continuum 2.}
\end{subfigure}
\hfill
\begin{subfigure}[b]{0.45\linewidth}
\includegraphics[width=\linewidth]{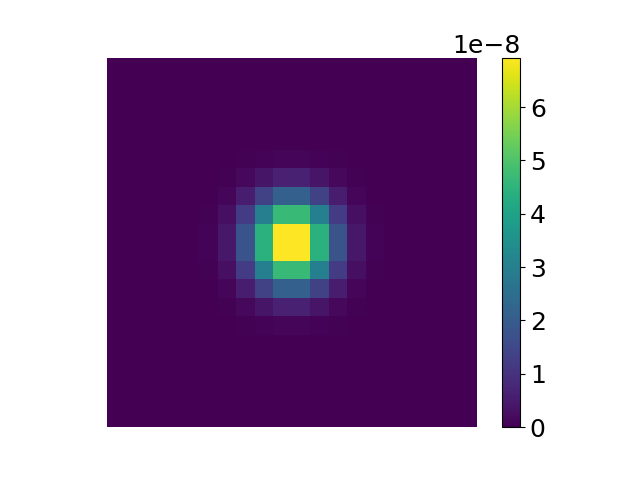}
\caption{Multiscale average in continuum 2.}
\end{subfigure}
\caption{Comparison of average solutions at final time in Example 1.}
\label{fig:ex1_uaUHa}
\end{figure}

To assess the fidelity of the multiscale approximation, Figure~\ref{fig:ex1_uaUHa} compares the continuum-wise average solutions computed from the fine-grid and multiscale simulations. The two solutions exhibit close agreement throughout the simulation, indicating that the multiscale model accurately captures the large-scale flow behavior governed by the underlying channel network.

\begin{table}[H]
\centering
\begin{tabular}{|c|c|c|}
\hline
\multicolumn{3}{|c|}{$H=1/10$} \\ \hline
$t$ & $e_2^{(1)}$ & $e_2^{(2)}$ \\ \hline
0 & 7.81e-02 & 8.42e-02 \\ \hline
1 & 4.24e-02 & 9.74e-02 \\ \hline
2 & 5.47e-02 & 6.92e-02 \\ \hline
3 & 5.75e-02 & 4.82e-02 \\ \hline
\end{tabular}
\hspace{0.5cm}
\begin{tabular}{|c|c|c|}
\hline
\multicolumn{3}{|c|}{$H=1/20$} \\ \hline
$t$ & $e_2^{(1)}$ & $e_2^{(2)}$ \\ \hline
0 & 1.97e-02 & 2.51e-02 \\ \hline
1 & 1.93e-02 & 2.83e-02 \\ \hline
2 & 5.05e-02 & 8.77e-03 \\ \hline
3 & 5.19e-03 & 4.86e-03 \\ \hline
\end{tabular}
\caption{Relative errors for different coarse grid sizes in Example 1.}
\label{tab:ers_ex1}
\end{table}

\subsection{Example 2}
This example considers a case where the second continuum $\kappa_2$ varies gradually throughout the domain. The numerical errors corresponding to different coarse grid sizes are reported in Table~\ref{tab:ers_ex2}. As the coarse mesh is refined, the error systematically decreases, indicating the expected convergence behavior. Moreover, the accuracy of the method is consistently preserved over the full simulation time.

\begin{figure}[H]
\centering
\begin{subfigure}[b]{0.45\linewidth}
\includegraphics[width=\linewidth]{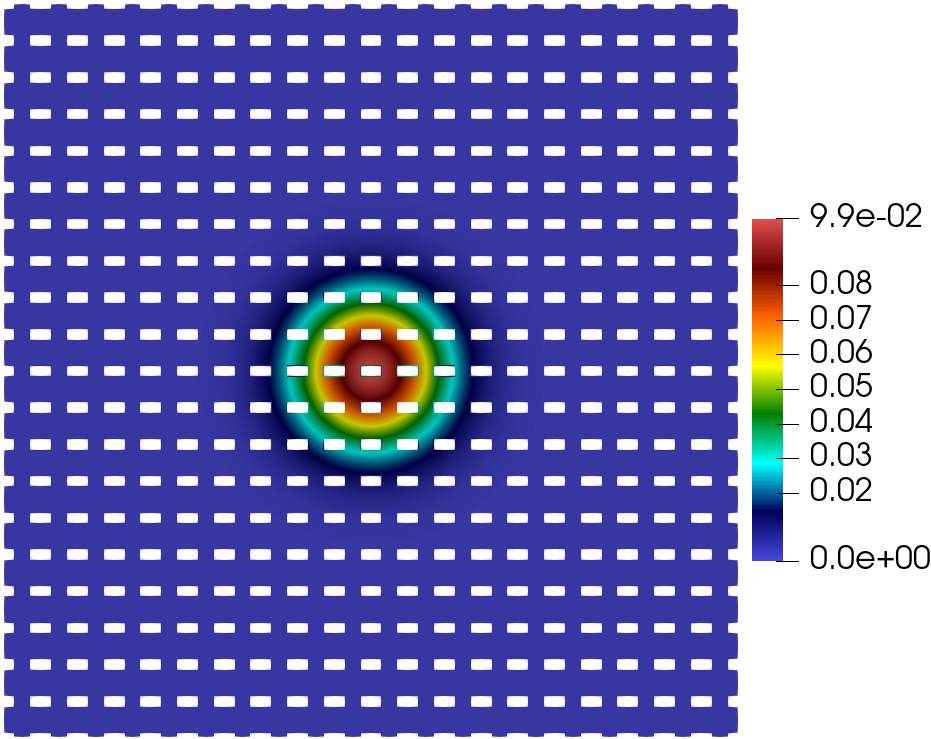}
\caption{$u_h^0$}
\end{subfigure}
\hfill
\begin{subfigure}[b]{0.45\linewidth}
\includegraphics[width=\linewidth]{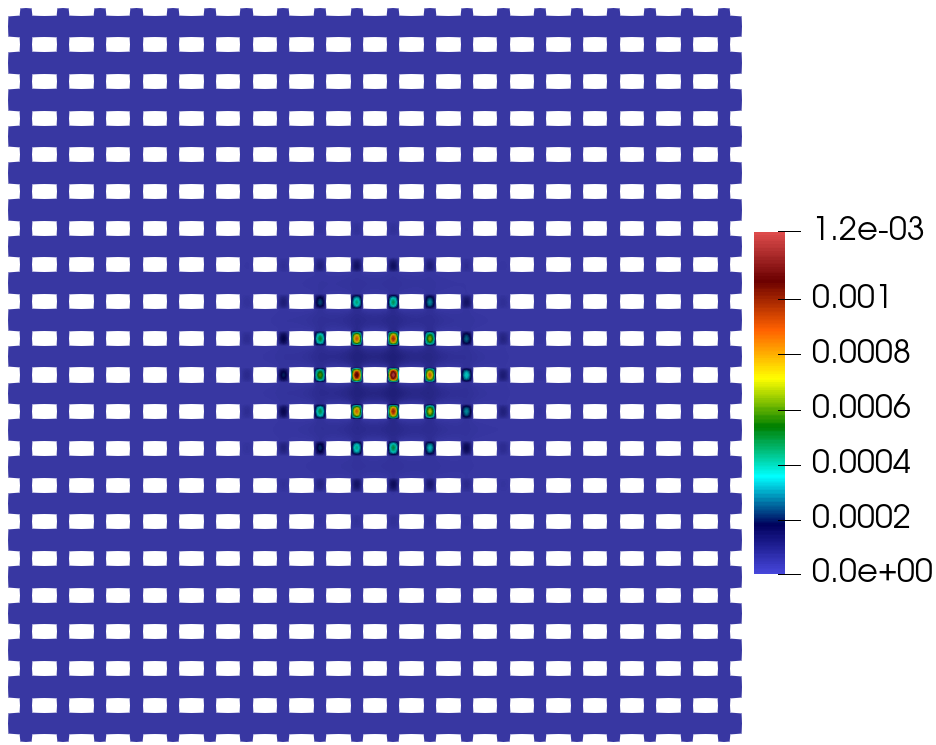}
\caption{$u_h^1$}
\end{subfigure}
\\[0.5cm]
\begin{subfigure}[b]{0.45\linewidth}
\includegraphics[width=\linewidth]{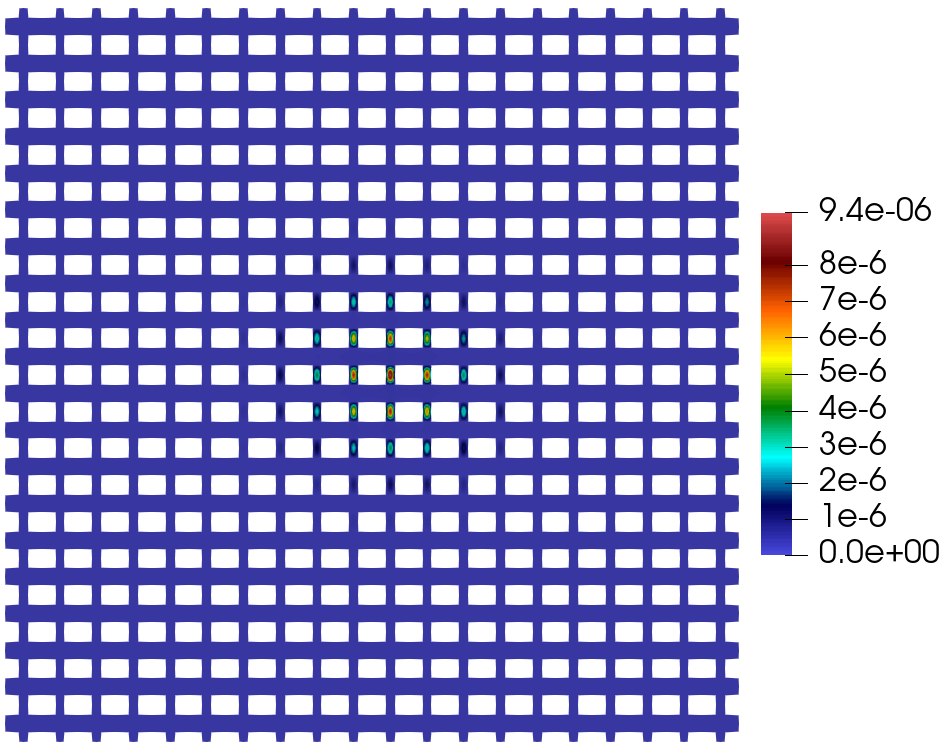}
\caption{$u_h^2$}
\end{subfigure}
\hfill
\begin{subfigure}[b]{0.45\linewidth}
\includegraphics[width=\linewidth]{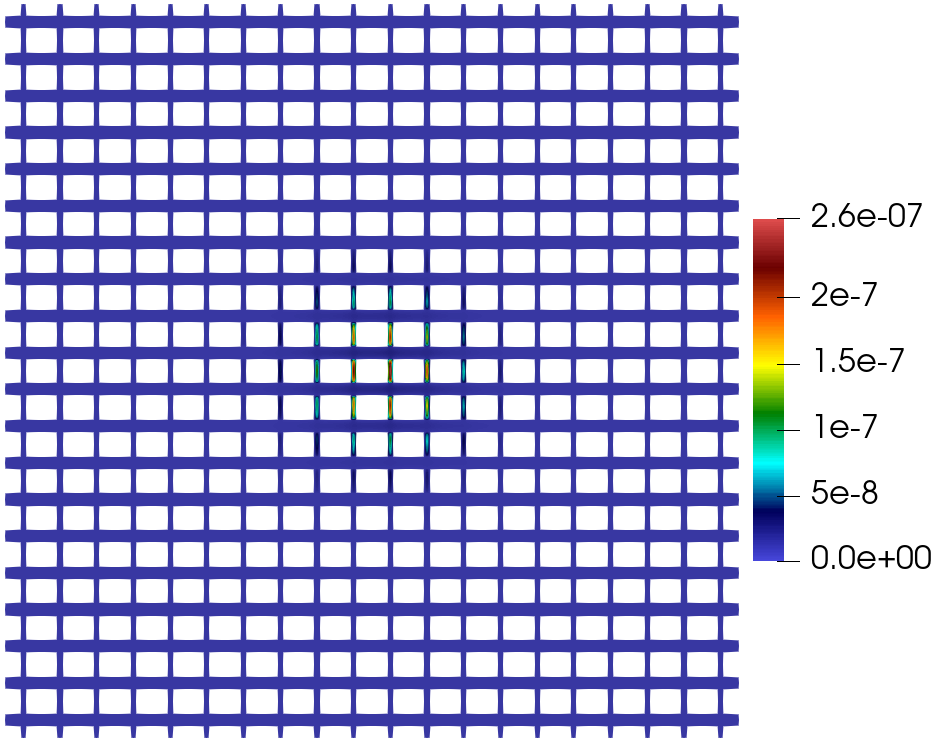}
\caption{$u_h^3$}
\end{subfigure}
\caption{Fine-grid solutions at different time steps in Example 2.}
\label{fig:ex2_uh}
\end{figure}

Figure~\ref{fig:ex2_uh} displays the fine-grid solutions at all time steps. The solution remains primarily localized within the narrow, low-permeability channels, where reduced conductivity restricts transport. In contrast, the high-permeability thick channels facilitate faster fluid movement, causing the solution to dissipate more quickly in those regions.

\begin{figure}[H]
\centering
\begin{subfigure}[b]{0.45\linewidth}
\includegraphics[width=\linewidth]{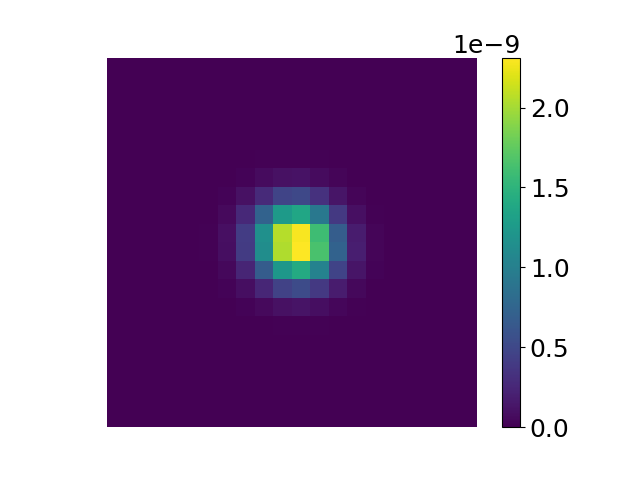}
\caption{Fine-grid average in continuum 1.}
\end{subfigure}
\hfill
\begin{subfigure}[b]{0.45\linewidth}
\includegraphics[width=\linewidth]{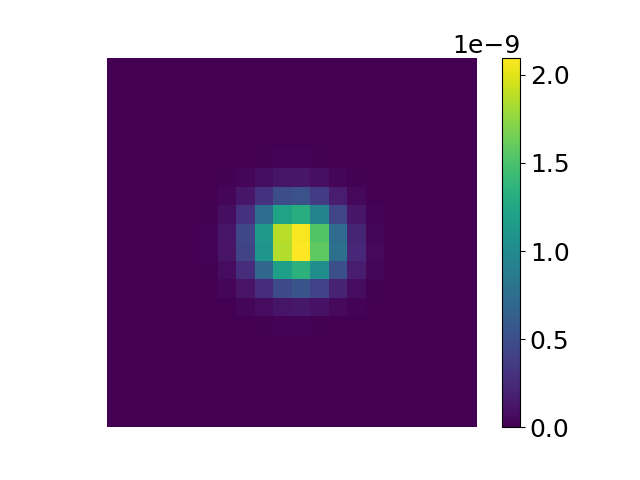}
\caption{Multiscale average in continuum 1.}
\end{subfigure}
\\[0.5cm]
\begin{subfigure}[b]{0.45\linewidth}
\includegraphics[width=\linewidth]{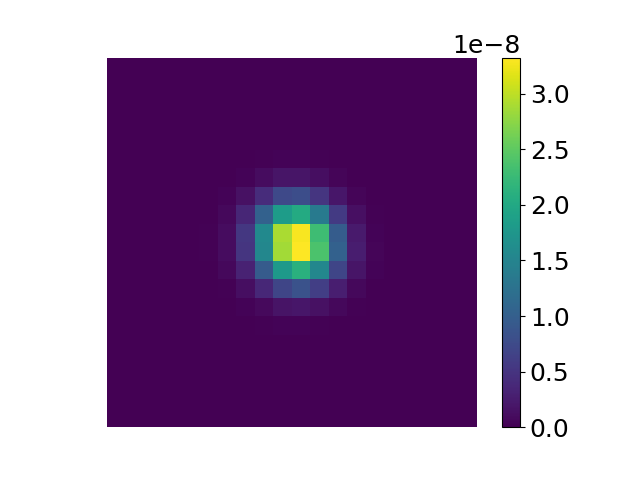}
\caption{Fine-grid average in continuum 2.}
\end{subfigure}
\hfill
\begin{subfigure}[b]{0.45\linewidth}
\includegraphics[width=\linewidth]{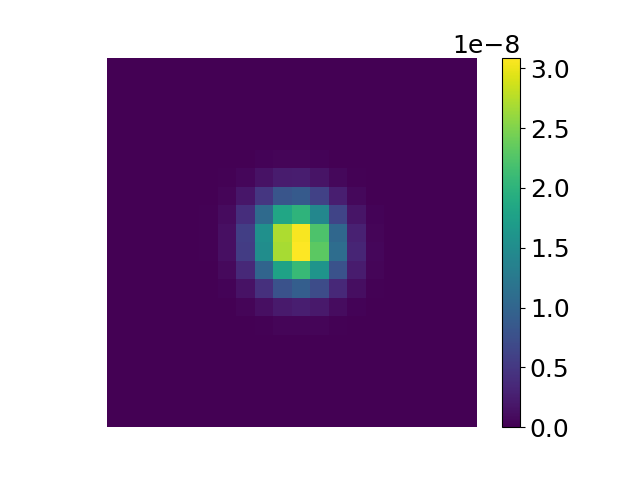}
\caption{Multiscale average in continuum 2.}
\end{subfigure}
\caption{Comparison of average solutions at final time in Example 2.}
\label{fig:ex2_uaUHa}
\end{figure}

To evaluate the performance of the multiscale method, Figure~\ref{fig:ex2_uaUHa} compares the averaged solutions obtained from the fine-grid and multiscale simulations. The two results show close similarity, with no substantial discrepancies observed. This demonstrates that the proposed method performs reliably under smoothly varying coefficient fields.

\begin{table}[H]
\centering
\begin{tabular}{|c|c|c|}
\hline
\multicolumn{3}{|c|}{$H=1/10$} \\ \hline
t & $e_2^{(1)}$ & $e_2^{(2)}$ \\ \hline
0 & 7.90e-02 & 8.58e-02 \\ \hline
1 & 4.45e-02 & 9.79e-02 \\ \hline
2 & 5.38e-02 & 7.03e-02 \\ \hline
3 & 5.85e-02 & 5.38e-02 \\ \hline
\end{tabular}
\begin{tabular}{|c|c|c|}
\hline
\multicolumn{3}{|c|}{$H=1/20$} \\ \hline
t & $e_2^{(1)}$ & $e_2^{(2)}$ \\ \hline
0 & 1.97e-02 & 2.51e-02 \\ \hline
1 & 1.97e-02 & 2.82e-02 \\ \hline
2 & 3.05e-02 & 8.58e-03 \\ \hline
3 & 5.25e-03 & 3.97e-03 \\ \hline
\end{tabular}
\caption{Errors for different coarse grid in example 2.}
\label{tab:ers_ex2}
\end{table}

\subsection{Example 3}

The numerical errors for the case under consideration are summarized in Table~\ref{tab:ers_ex3}. As shown, the error consistently decreases with the refinement of the coarse grid size \(H\). Notably, the accuracy of the solution is well-preserved across the entire simulation time interval.

\begin{figure}[H]
\centering
\begin{subfigure}[b]{0.45\linewidth}
\includegraphics[width=\linewidth]{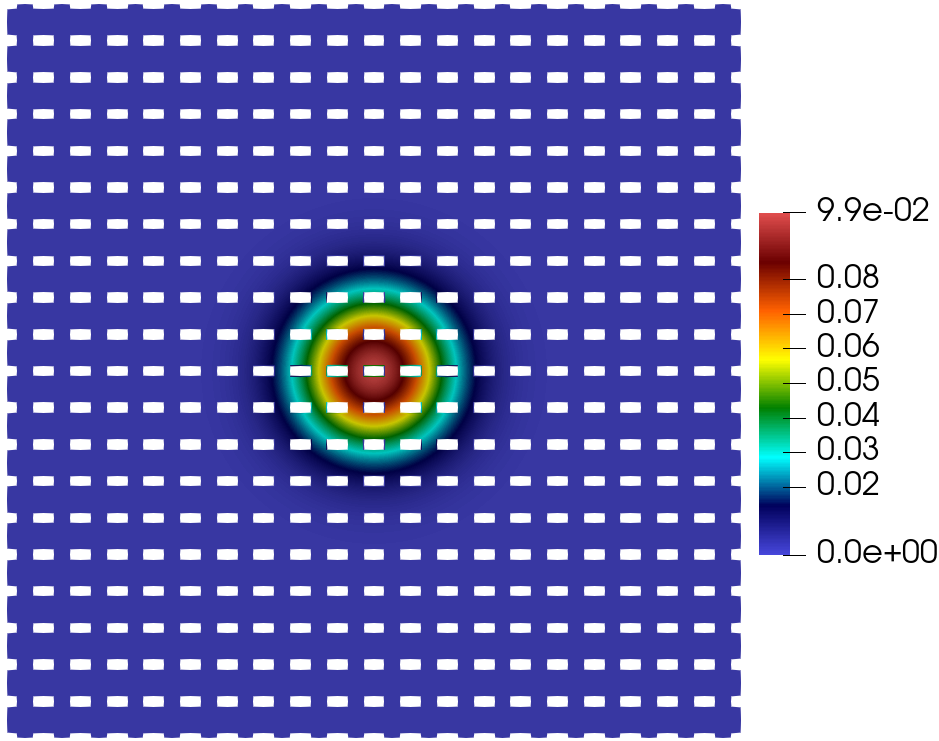}
\caption{$u_h^0$}
\end{subfigure}
\hfill
\begin{subfigure}[b]{0.45\linewidth}
\includegraphics[width=\linewidth]{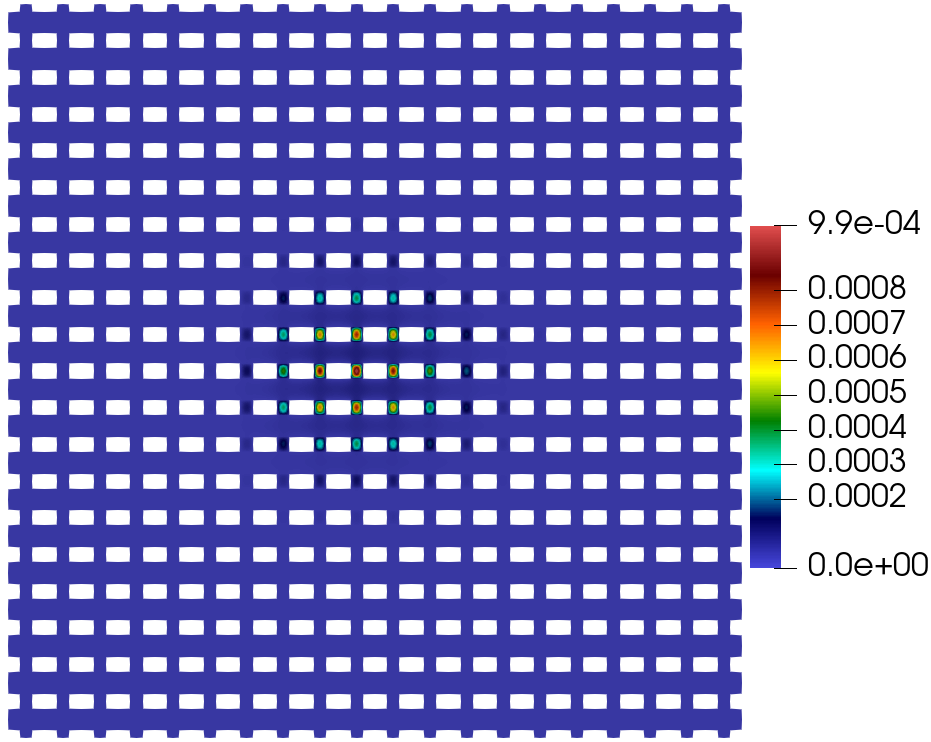}
\caption{$u_h^1$}
\end{subfigure}
\\[0.5cm]
\begin{subfigure}[b]{0.45\linewidth}
\includegraphics[width=\linewidth]{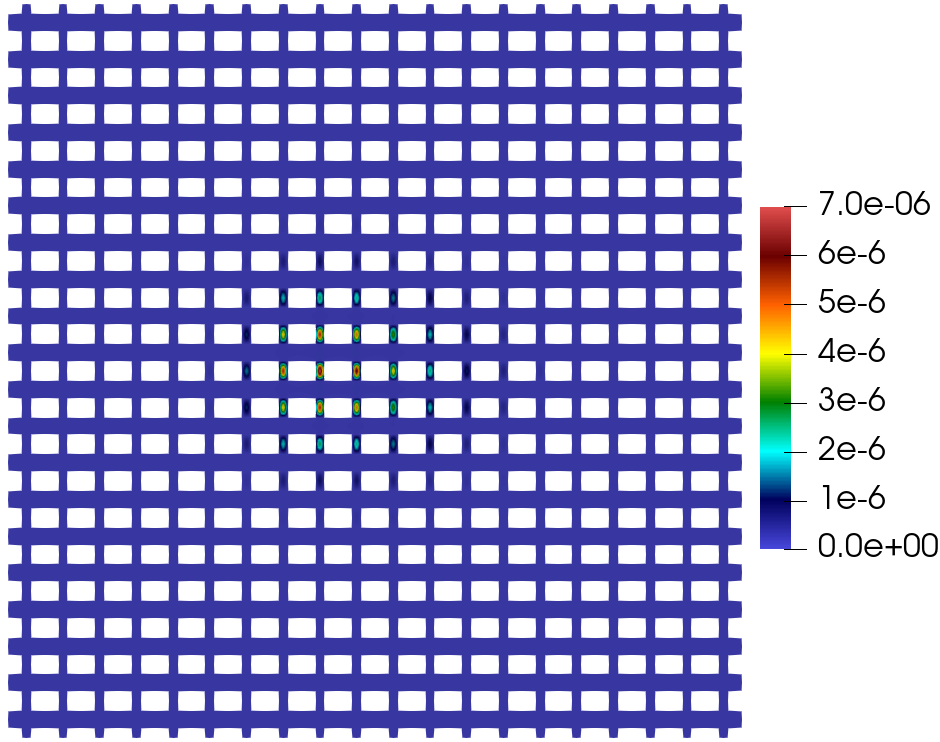}
\caption{$u_h^2$}
\end{subfigure}
\hfill
\begin{subfigure}[b]{0.45\linewidth}
\includegraphics[width=\linewidth]{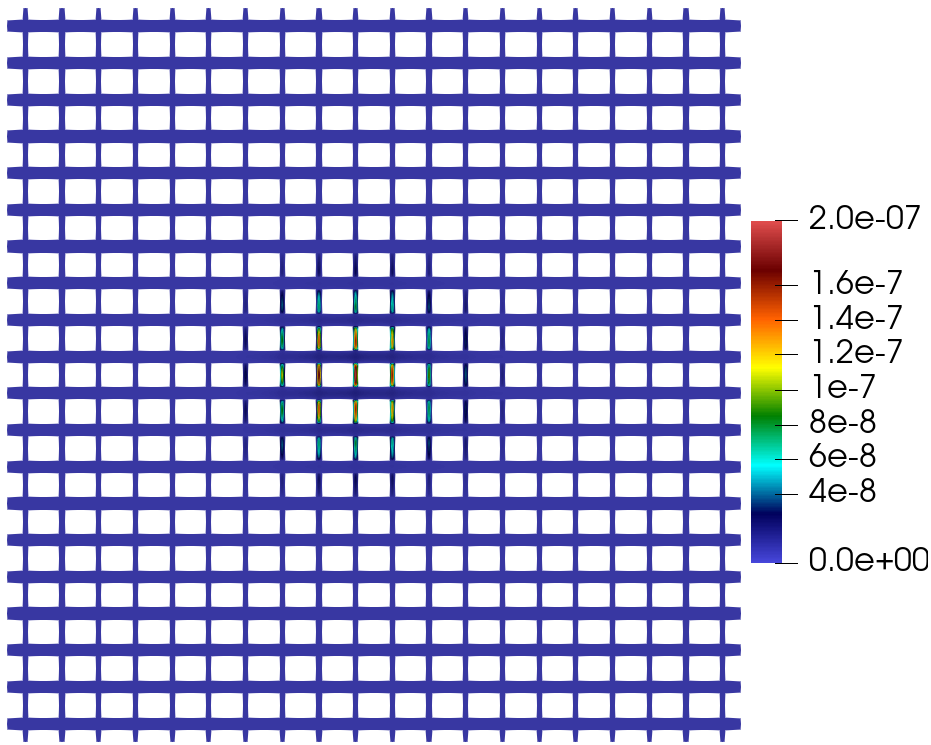}
\caption{$u_h^3$}
\end{subfigure}
\caption{Fine-grid solutions at different time steps in Example 1.}
\label{fig:ex3_uh}
\end{figure}

The fine-scale solutions at various time steps are illustrated in Figure~\ref{fig:ex3_uh}. It can be observed that the non-zero values are predominantly confined to the narrow channels, where the permeability is significantly lower than that of the broader channels. This behavior reflects the characteristic permeability contrast within the domain.

\begin{figure}[H]
\centering
\begin{subfigure}[b]{0.45\linewidth}
\includegraphics[width=\linewidth]{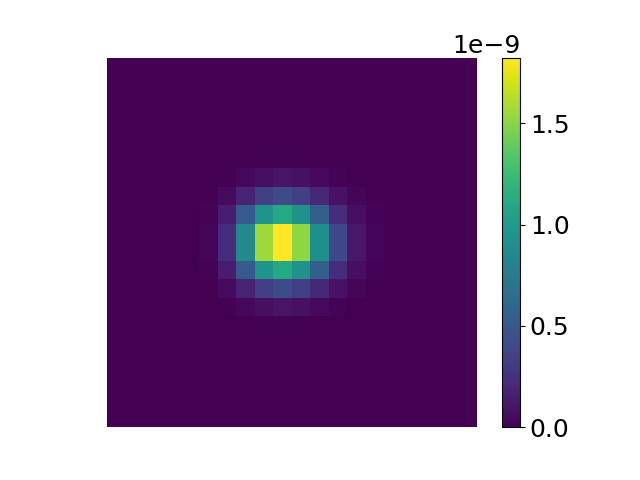}
\caption{Fine-grid average in continuum 1.}
\end{subfigure}
\hfill
\begin{subfigure}[b]{0.45\linewidth}
\includegraphics[width=\linewidth]{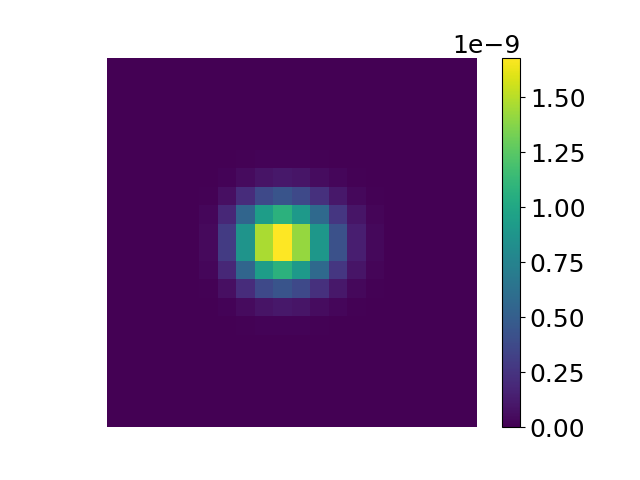}
\caption{Multiscale average in continuum 1.}
\end{subfigure}
\\[0.5cm]
\begin{subfigure}[b]{0.45\linewidth}
\includegraphics[width=\linewidth]{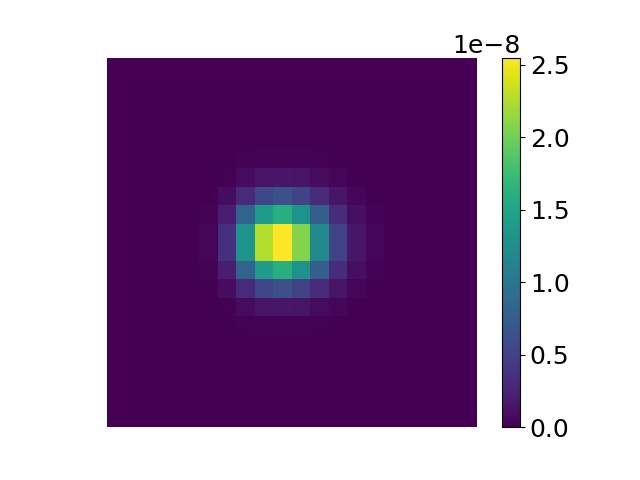}
\caption{Fine-grid average in continuum 2.}
\end{subfigure}
\hfill
\begin{subfigure}[b]{0.45\linewidth}
\includegraphics[width=\linewidth]{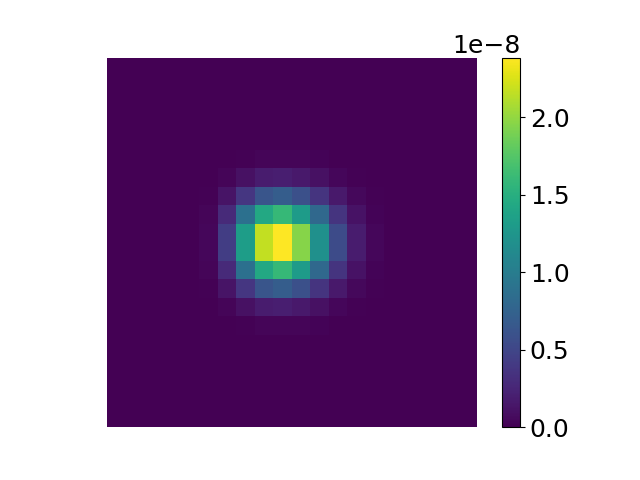}
\caption{Multiscale average in continuum 2.}
\end{subfigure}
\caption{Comparison of average solutions at final time in Example 1.}
\label{fig:ex3_uaUHa}
\end{figure}

Figure~\ref{fig:ex3_uaUHa} compares the averaged solutions obtained from both the fine grid and the multiscale method. The results exhibit negligible differences, further corroborating the efficiency and accuracy of the proposed multiscale approach.

\begin{table}[H]
\centering
\begin{tabular}{|c|c|c|}
\hline
\multicolumn{3}{|c|}{$H=1/10$} \\ \hline
t & $e_2^{(1)}$ & $e_2^{(2)}$ \\ \hline
0 & 8.55e-02 & 9.49e-02 \\ \hline
1 & 4.49e-02 & 8.67e-02 \\ \hline
2 & 5.36e-02 & 4.95e-02 \\ \hline
3 & 4.98e-02 & 4.54e-02 \\ \hline
\end{tabular}
\begin{tabular}{|c|c|c|}
\hline
\multicolumn{3}{|c|}{$H=1/20$} \\ \hline
t & $e_2^{(1)}$ & $e_2^{(2)}$ \\ \hline
0 & 1.99e-02 & 2.52e-02 \\ \hline
1 & 2.33e-02 & 2.58e-02 \\ \hline
2 & 3.50e-02 & 8.17e-03 \\ \hline
3 & 3.96e-03 & 3.72e-03 \\ \hline
\end{tabular}
\caption{Errors for different coarse grid in example 3.}
\label{tab:ers_ex3}
\end{table}

\section{Conclusions} \label{sec:conclusions}

In this work, we have developed a space-time upscaling method for problems involving dynamically evolving perforated domains. The proposed approach is based on multicontinuum homogenization and does not rely on scale separation assumptions. Distinct continua are identified according to the width of perforations, and corresponding cell problems with temporal derivatives are formulated over representative volume elements.
To illustrate the method, we consider parabolic problems, for which the proposed framework yields a coarse-grid model consisting of a coupled system of space-time equations defined on multiple continua.
Numerical experiments demonstrate that the method accurately captures the multiscale behavior of the solution, even in cases with complex spatial variations. The results confirm the robustness of the approach in handling dynamically changing microstructures.

\section*{Acknowledgement}
Wei Xie gratefully acknowledges the support from the China Scholarship Council (CSC, Project ID: 202308430231) for funding the research visit to Nanyang Technological University (NTU).
Viet Ha Hoang is supported by the Tier 2 grant  T2EP20123-0047 awarded by the Singapore Ministry of Education. 
Yin Yang is supported by the National Natural Science Foundation of China Project (No. 12261131501), the Project of Scientific Research Fund of the Hunan Provincial Science and Technology Department, China (No. 
2023GK2029, No. 2024JC1003, No. 2024JJ1008), and ``Algorithmic Research on Mathematical Common Fundamentals" Program for Science and Technology Innovative Research Team in Higher Educational Institutions of Hunan Province of China, and the 111 Project (No. D23017).

\bibliographystyle{abbrv}
\bibliography{references}

\end{document}